\magnification\magstep1
\hsize=16.8truecm
\vsize=25.4truecm
\baselineskip=12pt

\tolerance=1200
%
%
%
\long\def\title#1{\parindent 0pt
{\baselineskip 24pt\tolerance=10000
\vglue 3.9truecm
\noindent{#1}}\par}
\long\def\author#1{\vskip 12pt\begingroup{
\tolerance=10000\parindent=3truecm #1}\par\endgroup}

%
\countdef\sectno=11
\sectno = 0
\countdef\sbsectno=12
\sbsectno = 0
\countdef\ssbsectno=16
\ssbsectno=0
\def\section#1{\advance \sectno by 1 \sbsectno = 0
\sbsectno=0\ssbsectno=0\vskip 24pt
{\goodbreak
\noindent {\number\sectno.\ {#1}}}
\nobreak\vskip 12pt }
\def\subsection#1{\advance \sbsectno by 1
\ssbsectno=0
{\ifnum\count12=0\nobreak \else \medbreak \fi
\vskip 12 pt
\noindent {\number\sectno.\number\sbsectno. {\bf #1}}}
\nobreak\vskip 12pt }
\def\subsubsection#1{\advance \ssbsectno by 1
\medbreak\vskip 12pt{
\noindent{\number\sectno.\number\sbsectno.\number\ssbsectno.
\it #1.\ }}}
%
%
\countdef\figno=13
\figno = 0
\def\fig#1#2{\advance\figno by 1 \tolerance=10000
\setbox0=\hbox{\rm Fig. 00\ }
{\topinsert
\vskip #1 \smallskip\hangindent\wd0%
\noindent Fig.~\the\figno .\ #2
\vskip 12pt\endinsert}}
%
%
%

%
%
\def\REFERENCES
     {\countdef\refno=14
     \refno = 0\vskip 24pt
     {\centerline{\bf REFERENCES}\par\nobreak}
      \parindent=0pt
      \vskip 12pt\frenchspacing}
\def\ref#1#2#3
   {\advance\refno by 1
    \item{\rm[\number\refno] }\rm #1\ {\sl #2}
    \rm #3\par}
%
%
\def\frac#1#2{{#1 \over #2}}
%
%

%
%

%
%

%
%
\def\i{\ifmmode{\rm i}\else\char"10\fi}
%
%

%
\newfam\bofam
\font\tenbo=cmmib10   \textfont\bofam=\tenbo

\mathchardef\Omega="710A
\mathchardef\alpha="710B
\mathchardef\beta="710C
\mathchardef\gamma="710D
\mathchardef\delta="710E
\mathchardef\epsilon="710F
\mathchardef\rho="711A
\mathchardef\sigma="711B
\mathchardef\tau="711C
\mathchardef\upsilon="711D
\mathchardef\phi="711E
\mathchardef\chi="711F
\mathchardef\Gamma="7100
\mathchardef\Delta="7101
\mathchardef\Theta="7102
\mathchardef\Lambda="7103
\mathchardef\Xi="7104
\mathchardef\Pi="7105
\mathchardef\Sigma="7106
\mathchardef\Upsilon="7107
\mathchardef\Phi="7108
\mathchardef\Psi="7109
\mathchardef\zita="7110
\mathchardef\eta="7111
\mathchardef\theta="7112
\mathchardef\iota="7113
\mathchardef\kappa="7114
\mathchardef\lambda="7115
\mathchardef\mu="7116
\mathchardef\nu="7117
\mathchardef\xi="7118
\mathchardef\pi="7119
\mathchardef\psi="7120
\mathchardef\omega="7121
\mathchardef\varepsilon="7122
\mathchardef\vartheta="7123
\mathchardef\varpi="7124
\mathchardef\varrho="7125
\mathchardef\varsigma="7126
\mathchardef\varphi="7127

\baselineskip=18pt
\hsize=16truecm
\vsize=22truecm
\centerline{\bf ON THE BIFURCATION}
\centerline{\bf OF PERIODIC ORBITS.}
\bigskip\bigskip\bigskip
\centerline{\it J.P. Fran\c{c}oise}
\bigskip
\centerline{Universit\'e de Paris VI, GSIB, UFR de math\'ematiques,}
\centerline{ 175 Rue de Chevaleret,}
\centerline{75013 Paris, France}
{\vskip 12pt}

        In 1980, I visited IMPA and
J. Sotomayor asked me to "explain" Dulac's article on limit cycles.
Then Mauricio M. Peixoto explained me his interest in this question
in relation with his famous article on structural stability.
This was for me the
beginning of a fantastic experience with the world of differential
equations. The opportunity of Mauricio's $80^{th}$ birthday allows
me to express him all my gratitude  for his influence and to
warmly thank all his followers and close friends in particular
 Yvan Kupka, Charles C. Pugh, Jorge Sotomayor and Ren\'e Thom.

\bigskip
 \bigskip \bigskip\bigskip
\bigskip\bigskip
\noindent
{\bf INTRODUCTION}
{\vskip 12pt}

        Given a polynomial vector field of the plane:

 $$X=P(x,y){\partial}/{\partial x}+Q(x,y){\partial}/{\partial y},$$

 {\noindent}where $P$ and $Q$ are both polynomials of degree less than $d$, Hilbert's 16th problem, part B,
 calls for finding a bound (which depends only on $d$) to the number of limit cycles (isolated
periodic solutions) of $X$ . Despite the many contributions to this question, the problem still
remains open nowadays even in the case $d=2$.
        One of the most spectacular contribution to this field is due to N.N. Bautin ([1]).
These notes aim at providing generalizations of Bautin's theorem which involve several
new ideas introduced in recent years.
        The creation and developments of the theory of non linear oscillations is due to the school of
A.A.Andronov. The famous books " Theory of oscillations " by Andronov, Khaikin, Witt" followed by
"Qualitative theory of second order dynamical systems" and "Theory of bifurcations of dynamical
systems in the plane" run several editions in various languages. As one of the most prominent
followers of Andronov, N.N. Bautin wrote in 1939 a masterpiece contribution to the classical
center-focus problem posed by Poincar\'e.
{\indent}In his lecture given on the occasion of Arnol'd's 60-th birthday at the Fields Institute in
Toronto, June 1997, ([13], [14]) S. Smale listed 18 problems chosen with these criteria:
simple statement, also preferably mathematically precise; personal acquaintance
with the problem
; a belief that the question, its solution, partial results, or even attempts at its solution are
likely to have great importance for Mathematics and its developments in the next century.
        Besides the Riemann Hypothesis, one was on Hilbert's list; the Hilbert's 16-th problem on limit
cycles.
        The author hopes that the specialists will be interested to have here collected techniques which
appeared in a scattered way (from 1995 to 2000) and have been revisited in a more systematic manner.
For non-specialists, this survey may contribute to introduce this vivid subject which involves
so different methods such as projections of analytic sets, division by an ideal, complex
foliations, Riemann surfaces, algebraic invariant theory and normal forms, and which belongs to the
tradition of the qualitative theory of differential equations.
Chapters I-II review results which have essentially already been published
jointly with Yosi Yomdin ([8], [9]).
  Chapter III is perhaps the most original
and develops an extension of Bautin's approach to any dimension in the framework of the
 theory of pertubations of a Nambu dynamics ([12]).

\vfill\eject

\vskip 12pt
\noindent
{\bf I- CLASSICAL BAUTIN's APPROACH}
\vskip 12pt

        The first paragraph is devoted to the presentation of Bautin approach ([1]). We give the definitions of the center set, of
the Bautin ideal and of the Bautin index. We show that the center set is a real algebraic
variety. We select particular generators of the Bautin ideal which are adapted to the first
return map in a sense to be precised. We give an estimate of the domain of convergence of
the return mapping and we show accordingly to ([2], [3], [4]) that it is an $A_{0}$-Series.
We conclude this part with Bautin's argument based on the classical Rolle property to produce
a bound to the number of real zeros of the return mapping.
        In this survey, we will mostly be concerned with homogeneous nonlinearities.
 In this paragraph, it is not necessary to assume that the perturbation is homogeneous: all
definitions and theorems can be easily modified without any fundamental changes.
Consideration of this particular case allows simplification of
the notations and also it captures the essential properties of the general
case.

        \vskip 12pt
\noindent
{\bf I.1 Construction of the return mapping by Bautin's method}
\vskip 12pt

        We consider a polynomial plane vector field $X$ of type:

$$X=x{\partial}/{\partial y}-y{\partial}/{\partial x}+
 \sum_{i,j/ i+j=d}
 [ a_{i,j} x^{i}y^{j}{\partial}/{\partial x} +
  b_{i,j} x^{i}y^{j}{\partial}/{\partial y}].
  \eqno(1.1)$$
and the associated flow, solution of the system:

$$\dot{x}= -y +\sum _{i,j/ i+j=d} a_{i,j} x^{i}y^{j}= -y+P(x,y), \eqno(1.2a)$$
$$\dot{y}= x + \sum_{i,j/ i+j=d} b_{i,j} x^{i}y^{j}= -y+Q(x,y). \eqno(1.2b)$$

        The parameters of the vector field $(a,b)$ can take any real values and thus $(a,b)$ should be
considered as a point in the vector space $R^{2(d+1)}$. We first recall Bautin's approach to
find the return mapping of (1.1) near the origin.
        Write (1.2) in polar coordinates $(r,\theta)$:

 $$ x = r{\rm cos}(\theta),  y = r{\rm sin}(\theta). \eqno(1.3)$$

        This leads to:

 $$ 2r{\dot r} = 2(x{\dot x}+y{\dot y}), r{\dot r}=xP+yQ= r^{d+1}A(\theta), \eqno(1.4)$$

 $$ {\dot \theta} = (x{\dot y}-y{\dot x})/(x^{2}+y^{2}) = 1+ r^{d-1}B(\theta), \eqno(1.5)$$

 where $A(\theta)$ and $B(\theta)$ are two trigonometric polynomials (in ${\rm cos}(\theta),
 {\rm sin}(\theta)$) linear in the parameters $(a,b)$.

        This yields:

  $${dr}/{d\theta} = r^{d}A(\theta)/[1+r^{d-1}B(\theta)], \eqno(1.6)$$

  and thus

  $${dr}/{d\theta} = \sum_{k=0}^{\infty} (-1)^{k}r^{k(d-1)+d}A(\theta)B(\theta)^{k}. \eqno(1.7)$$

  This equation (1.7) may be rewritten:

  $${dr}/{d\theta} = \sum_{k\geq d} r^{k}R_{k}(\theta), \eqno(1.8)$$

  where the coefficients $R_{k}(\theta)$ are trigonometric polynomials in $({\rm cos}(\theta),
 {\rm sin}(\theta))$ and polynomials in the parameters $(a,b)$. To simplify the notations, the
dependence on the parameters $(a,b)$ is not made explicit.
        Bautin looks for a solution of (1.8) $r = r(\theta)$ so that $r(0) = r_{0}$, given as an expansion:

 $$r = r_{0} + v_{2}(\theta)r_{0}^{2}+ ... +  v_{k}(\theta)r_{0}^{k}+ ...\eqno(1.9)$$

 Comparison between (1.8) and (1.9) yields:

 $$ v_{2}(\theta) = ... =  v_{d-1}(\theta) = 0, \eqno(1.10)$$
 $$ v_{d}'(\theta) = R_{d}(\theta), \eqno(1.11)$$
 $$ v_{k}'(\theta) = \sum^{k}_{i=2} B_{ik}[v_{d}(\theta), ... , v_{k-1}(\theta)]R_{i}(\theta),
 k \geq d+1. \eqno(1.12)$$
  The polynomial $B_{ik}[a_{d}, ... , a_{k-1}]$ displays integer coefficients and is the coefficient
of $X^{k-i}$ in $(X+a_{d}X^{d}+...+a_{p}X^{p}+...)^{i}$.
        The relation (1.12) allows to determine inductively the functions $v_{k}(\theta)$:

 $$ v_{k}(\theta) = \int^{\theta}_{0} [\sum^{k}_{i=2} B_{ik}[v_{d}(\phi), ... ,
v_{k-1}(\phi)]R_{i}(\phi)]d\phi. \eqno(1.13)$$

        Two facts can be easily derived from this construction:

 i) $v_{k}(\theta)$ is polynomial in $\theta$ (of degree less than $k$) and in $(sin(\theta),
cos(\theta))$.

ii) $v_{k}(\theta)$ is polynomial in the parameters $(a,b)$ of the vector field. Thus in particular
the coefficients $v_{k}(2\pi)$ of the return mapping are polynomials in $(a,b)$.
\vskip 12pt
{\bf Definition I.1.1}
\vskip 12pt

        The vector field $X$ displays a center at the origin $(0\in R^{2}$)
 (or is said to be a center) if and only if $r(2\pi)=r(0)$
  for all values of $r(0)$ close to $0$.

\vskip 12pt
\noindent
{\bf I.2 Convergence of the Taylor series, majorant series and $A_{0}$-series}
\vskip 12pt

        Let $ f_{\lambda}(x)= \sum a_{k}(\lambda)x^{k}$ be an analytic series in $x$ with
polynomial coefficients in the parameters $\lambda= (\lambda_{0},...,\lambda_{D})$. Denote
$\mid a_{k}\mid$ (norm of the polynomial $a_{k}(\lambda)$) as the sum of the absolute
value of the coefficients.
\vskip 12pt
{\bf Definition I.2.1}
\vskip 12pt

        The series $f_{\lambda}$ is called an $A_{0}$-series if the following two conditions are satisfied:

 There are positive constants $K_{1}, K_{2},K_{3},K_{4}$ such that:

 1- $deg(a_{k})\leq K_{1}k+K_{2},$

 2- $\mid a_{k}\mid \leq K_{3}K_{4}^{k}$.
 \vskip 12pt

 $A_{0}$-series form a subring of the ring of formal power series in $x$ with polynomial
coefficients in $\lambda$. All the usual analytic operations, like substitution to a given analytic
function, composition, inversion, etc... transform $A_{0}$-series into themselves. (The proof is
rather straightforward and is not provided here). The $A_{0}$-series have been precisely introduced (in the
subject) by M. Briskin and Y. Yomdin ([4]). The proof of the following easy lemma is also left to the reader:
\vskip 12pt
{\bf Lemma I.2.2}
\vskip 12pt

        An $A_{0}$-series $f_{\lambda}(x)$ converges in the disc $D(0,R)$ of radius $R=
[1/(K_{4}(1+\mid\lambda\mid))^{K_{1}}]$.
\vskip 12pt
        In the following, we also denote by  $f_{\lambda}$ the complex analytic function defined for all
$\lambda \in C^{D}$ on the disc $D(0,R)$ by the $A_{0}$-series.
\vskip 12pt

{\bf Proposition I.2.3}
\vskip 12pt

        For all $\theta$, the series $\sum_{k\geq d} x^{k}R_{k}(\theta),$ is an $A_{0}$-series
 with $K_{1}= 1/(d-1)$, $K_{2}= -1/(d-1)$, $K_{3}= 1/[2(d+1)]^{1/d-1}$, $K_{4}= [2(d+1)]^{1/d-1}$.

\vskip 12pt
{\bf Proof.}
\vskip 12pt

        This is indeed a simple consequence of Bautin's method. For all $\theta$, the norms of the
polynomials $A(\theta)$ and $B(\theta)$ (seen as polynomials in
$(a,b)$),

 $$A(\theta)= (xP+yQ)[cos(\theta), sin(\theta)], \eqno(1.14 a)$$

$$B(\theta)= (xQ-yP)[cos(\theta), sin(\theta)], \eqno(1.14 b)$$

{\noindent}are bounded by:

$$\mid A(\theta)\mid \leq 2(d+1), \mid B(\theta)\mid \leq 2(d+1). \eqno(1.15)$$
Write:

$${dr}/{d\theta} = \sum_{k\geq d} r^{k}R_{k}(\theta) = \sum_{j\geq 0}
(-1)^{j}r^{d+j(d-1)}A(\theta)B(\theta)^{j}. \eqno(1.16)$$

Denote $k=d+j(d-1)$, then this yields:

$$ deg[R_{k}(\theta)] \leq 1+j \leq 1+ [(k-d)/(d-1)] \leq [(k-1)/(d-1)]. \eqno(1.17)$$

Furthermore, the norm of $R_{k}(\theta)$ as polynomial in the parameters $(a,b)$ is estimated as
follows:

$$\mid R_{k}(\theta)\mid = \mid A(\theta)B(\theta)^{j}\mid \leq \mid A(\theta)\mid\mid
B(\theta)\mid^{j} \leq
[2(d+1)]^{[(k-1)/(d-1)]}.
\eqno(1.18)$$
\vskip 12pt

{\bf Theorem I.2.4}
\vskip 12pt

                For all $\theta$, the series $\sum_{k\geq d} x^{k}v_{k}(\theta),$ is an $A_{0}$-series
 with $K'_{1}= 1/(d-1)$, $K'_{2}= 0$, $K'_{3}=
[2\pi K_{3}/4K_{4}^{2}][K_{4}-2C+((K_{4}-2C)^{2}-K_{4}^{2}]^{2}$, $1/K'_{4}=
|K_{4}-2C+((K_{4}-2C)^{2}-K_{4}^{2}|
 /[2K_{4}^{2}]$, $C= K_{4}+2\pi K_{3}K_{4}^{2}$.
\vskip 12pt
{\bf Proof.}
\vskip 12pt
  First observe that $deg(v_{d}(\theta)) = deg(R_{d}(\theta)) = 1$ (The degree as polynomial in the
parameters). Thus we have $deg(v_{d}(\theta)) \leq d/(d-1)$. Assume inductively that:

$$deg(v_{j}(\theta)) \leq j/(d-1) for j=d,...,k-1. \eqno(1.19)$$

        The recurrency relation (1.13) displays:

 $${\rm deg}[v_{k}(\theta)]\leq max [ {\rm deg}(B_{ik}[v_{d}(\phi)), ... ,
{\rm deg}(v_{k-1}(\phi)]R_{i}(\phi))], i=2,...,k. \eqno(1.20)$$

        The equation (1.17) and the induction assumption yield:

$$deg [B_{ik}[v_{d}(\phi), ... ,v_{k-1}(\phi)]R_{i}(\phi)] \leq K_{1}(k-i)+K_{1}i+K_{2} \leq K_{1}k.
\eqno(1.21)$$

        This shows the first part of the theorem on the bound of the degrees of the coefficients
$v_{k}(\theta)$.
        For the second part of the proof related to the bound on the norms of the coefficients, we use
standard methods and notations of majorant series.
\vskip 12pt

{\bf Definition I.2.5}
\vskip 12pt

        The formal series $\Psi(x)= \sum_{k \geq 1}\Psi_{k}x^{k}$ with positive coefficients dominates the
formal series $\Phi(x)= \sum_{k \geq 1}\Phi_{k}x^{k}$ with positive coefficients if and only if for
all coefficients $\Phi_{k} \leq \Psi_{k}, k \geq 1$.
\vskip 12pt
        The series $x+\sum_{k \geq 2} {\mid v_{k}(\theta)\mid x^{k}}$ is dominated by the series
 $x+\sum_{k \geq 2} {\Psi_{k}x^{k}}$ so that:

 $$ \psi_{k} = 2\pi\sum_{i=d}^{k} B_{ik}(\psi_{d},...,\psi_{k-1})\mid R_{i}\mid.\eqno(1.22)$$

 Denote
 $$\bar{R}(x)= \sum_{i \geq 2} K_{3}K_{4}^{i}x^{i} = K_{3}(K_{4}x)^{2}/[1-K_{4}x], \eqno(1.23)$$

 then the proposition I.2.3 yields:

 $$\sum_{i \geq d}\mid R_{i}\mid x^{i} {\rm is} {\rm dominated} {\rm by} \bar{R}(x). \eqno(1.24)$$

        The series $\Psi(x) = x+\sum_{k \geq 2} {\Psi_{k}x^{k}}$ is then dominated by the solution
  $\bar{\Psi}(x)$ solution of the equation:

 $$ \bar{\Psi}(x) - x = 2\pi \bar{R}[\bar{\Psi}(x)]=
2\pi K_{3}[K_{4}\bar{\Psi}(x)]^{2}/[1-K_{4}\bar{\Psi}(x)]. \eqno(1.25)$$

At this point we have obtained that $\bar{\Psi}(x)$ is a solution of an algebraic equation of
degree two. Estimates of the constant $K_{4}'$ is then obtained by the distance to the first zero
of the discriminant and constant $K'_{3}$ is then adjusted from the first term of the development.
Easy computations give then the proof of the theorem.
\vskip 12pt
\noindent

 {\bf I.3 Definitions of the center set, of the Bautin ideal and of the Bautin index. Definition of
$\phi$-equivalence of series.}
\vskip 12pt

{\bf Definition I.3.1}
\vskip 12pt

        The center set $C$ is the set of values of parameters $(a,b)$ so that the corresponding
vector field $X$ has a center at the origin.
\vskip 12pt

{\bf Proposition I.3.2}
\vskip 12pt

        The center set $C$ is a real algebraic manifold.
\vskip 12pt
{\bf Proof}
\vskip 12pt

 The center set $C$ is given as the zero set of the coefficients $v_{k}(2\pi)$ which are
polynomials in the parameters $(a,b)$.
\vskip 12pt

        From now on, it is appropriated to change of notations and denote $L_{k}(a,b)=v_{k}(2\pi)$ the
coefficients of the return mapping to emphasize their dependence in terms of the parameters $(a,b)$.
We denote:

$$r \mapsto L(r) = r+ L_{d}(a,b)r^{d}+...+L_{k}(a,b)r^{k}+... \eqno(1.26)$$
the return mapping defined in (1.9) for $\theta=2\pi$.
\vskip 12pt
{\bf Definition I.3.3}
\vskip 12pt

        The Bautin ideal is the ideal generated in the ring $R[a,b]$ by the coefficients $L_{k}(a,b)$.

\vskip 12pt
{\bf Definition I.3.4}
\vskip 12pt

        The Bautin index is the first integer $k_{0}$ so that the polynomials
$L_{d}(a,b),...,L_{k_{0}}(a,b)$ generate the Bautin ideal.
\vskip 12pt

        The reader should be careful about the fact that the Bautin index does not depend only on the
Bautin ideal. In its definition, we cannot substitute to the collection of the coefficients of the
return mapping another system of generators of the Bautin ideal.
        Note that the existence of the Bautin index just follows from the fact that the ring $R[a,b]$ is
Noetherian.

                The local Hilbert's 16th problem is to find a bound depending only on $d$ to the number of
isolated periodic orbits of $X$ in a neighborhood of the origin. Isolated periodic orbits of $X$
defined in the domain of definition of the return mapping correspond exactly to the isolated
solutions of the equation:

$$L(r)-r=0.\eqno(1.27)$$

        We gradually change into notations more pertinent to the general algebraic geometry setting. Let
 $\Phi(x,\lambda)$ be an analytic series in $x$ with polynomial coefficients in the parameter
$\lambda=(\lambda_{1},...,\lambda_{D})$.
\vskip 12pt
{\bf Definition I.3.5}
\vskip 12pt

        The Bautin ideal of the series $\Phi(x,\lambda)$ is the ideal generated in the ring $R[\lambda]$
 by the coefficients $\Phi_{k}(x,\lambda)$. The center set of the series $\Phi(x,\lambda)$ is the
zero set of its Bautin ideal. The Bautin index $d$ of the series $\Phi(x,\lambda)$ is the minimal
integer $d$ such that the coefficients $\Phi_{1}(x,\lambda),...\Phi_{d}(x,\lambda)$ generate the
Bautin ideal of the series $\Phi(x,\lambda)$.
\vskip 12pt

{\bf Definition I.3.6}
\vskip 12pt

        Two series $\Phi(x, \lambda)$ and  $\Psi(x, \lambda)$ with polynomial coefficients in the
parameters $\lambda$ are said to be $\phi$-equivalent if for all integers $k \geq 1$, the polynomial
$\Phi_{k}(\lambda)-\Psi_{k}(\lambda)$ belongs to the ideal generated by $(\Phi_{1}(\lambda),...,
\Phi_{k-1}(\lambda))$.
\vskip 12pt

        It can be easily checked that this defines an equivalence relation. Also the definition yields the
following:
\vskip 12pt
{\bf Lemma I.3.7}
\vskip 12pt

        Two series $\Phi(x, \lambda)$ and  $\Psi(x, \lambda)$ with polynomial coefficients in the
parameters $\lambda$ which are $\phi$-equivalent have the same Bautin index.

\vskip 12pt
\noindent
{\bf I.4 The classical argument to show that the number of real zeros is bounded}
\vskip 12pt

        We recall here Bautin's argument to produce a bound to the number of real zeros.

 Let $\Phi(x,\lambda)$ be an analytic series in $x$ with polynomial coefficients in the parameters
 $\lambda$ and with Bautin index $d$. Assume that $\Phi(0,\lambda) =0$ and $\Phi(x,0) =0$.
\vskip 12pt
 {\bf Theorem I.4.1}
 \vskip 12pt

        There is a ball $B\in R^{D}$, centered at 0, in the space $R^{D}$ and an interval $I$ containing
$0$ such that for all $\lambda \in B$, the number of zeros of $\Phi(x,\lambda)$ contained in $I$ is
less than or equal to $d$.
\vskip 12pt
{\bf Proof.}
\vskip 12pt

        Using the definition of the Bautin index $d$, we write:

 $$\Phi(x,\lambda)= \sum^{d}_{i=1}{\Phi_{i}(\lambda)[1+\Psi_{i}(x,\lambda)]x^{i}}, \eqno(1.28)$$

 with $\Psi_{i}(x,\lambda)$ analytic in $x$, with polynomial coefficients in $\lambda$ such that
 $\Psi_{i}(0,\lambda)=\Psi_{i}(x,0)=0$. Assume that $B$ and $I$ are small enough so that (for
instance)  $\mid \Psi_{i}(x,\lambda)\mid \leq 1/2$ on $I\times B$. Then divide $\Phi(x,\lambda)$ by
$[1+\Psi_{d}(x,\lambda)]$ and rewrite:

$$[\Phi(x,\lambda)]/[1+\Psi_{d}(x,\lambda)] = \Phi_{1}(\lambda)+
\Phi_{2}(\lambda)[1+\Psi'_{i}(x,\lambda)]x^{2}+... \eqno(1.29)$$

        Then from Rolle's lemma follows that the number of zeros of $\Phi(x,\lambda)$ is less than 1+
number of zeros of the derivative $[\Phi(x,\lambda)]/[1+\Psi_{d}(x,\lambda)]'$. Write then this
derivative as

$$ \Phi_{2}(\lambda)[1+\Psi^{(2)}_{i}(x,\lambda)]x+...
\Phi_{d-1}(\lambda)[1+\Psi^{(2)}_{i}(x,\lambda)]x^{d-1}. \eqno(1.30)$$

        Then repeat the process (sometimes referred to as the division-derivation algorithm). We obtain the
result by an easy induction.
\vskip 12pt
        It is clear that one cannot go far with this proof if we want an explicit control of the size of
the domain (either of $I$ or of $B$). For instance, we know that the coefficients
$\Phi_{i}(\lambda)$ are combinations of the $d$ first coefficients:

$$\Phi_{i}(\lambda) = P_{(i,1)}(\lambda)\Phi_{1}(\lambda)+...+P_{(i,d)}(\lambda)\Phi_{d}(\lambda),
\eqno(1.31)$$

{\noindent}but how to control the size of the "quotients" $P_{(i,k)}(\lambda)$ which enter in the construction
of the $\Psi_{i}(x,\lambda)$? We will see how to bypass these difficulties in the next chapter
using complex analytic methods such as the polynomial Hironaka division theorem.
\vskip 12pt
\noindent

{\bf I.5 Formulation of the problem in terms of projection of analytic sets.}
\vskip 12pt

        Let $\Phi:R\times R^{D} \mapsto R$ be an analytic series with polynomial coefficients:

 $$\Phi(x,\lambda)= x+\Phi_{2}(\lambda)x^{2}+...+\Phi_{k}(\lambda)x^{k}+... \eqno(1.32)$$

        We consider the subset $\Sigma\subset R\times R^{D}$ defined as the zero-set of $\Phi(x,\lambda)-x$:

 $$ \Sigma = [ (x,\lambda)/ \Phi(x,\lambda)-x = 0]. \eqno(1.33)$$

        Let $\pi :\Sigma \mapsto R^{D}$ be the restriction to $\Sigma$ of the natural projection
  $\pi: R\times R^{D} \mapsto R^{D}$. The center set associated to the analytic series
$\Phi(x,\lambda)-x$ defined in I.3.5 is the set $C\subset R^{D}$ of parameters $\lambda$ such that the
fiber of the projection $\pi^{-1}(\lambda)$ is contained in the set $\Sigma$.
        The theorem I.4.1 displays the following geometric interpretation:
\vskip 12pt
{\bf Theorem I.5.1}
\vskip 12pt
        There is a neighborhood $I\times B$ of $(0,0)$ in $R\times R^{D}$ such that for all points $\lambda$ of
  $B$ the number of isolated points of the fibers $\pi^{-1}(\lambda)$ restricted to $\Sigma$
  is less than the Bautin index $d$ of the analytic series $\Phi(x,\lambda)-x$.
\vskip 12pt
\noindent
  {\bf I.6 A $S^{1}$-action on the space of parameters and a system of $S^{1}$-invariant
polynomials which generate the Bautin ideal.}
\vskip 12pt

        We can change the coordinates in the plane $(x,y) \mapsto (x',y')$ under the action of the
rotation group. This leaves invariant the linear part of the vector field $X$ and thus it induces
an action of the rotation group $S^{1}$ on the space of parameters $(a,b)$. This linear action
induces an action on the ring of polynomials $R[a,b]$ on the space of parameters. A polynomial is
said to be invariant if it is fixed under this $S^{1}$-action. The center set is obviously
invariant under the $S^{1}$-action because the fact, that the vector field $X$ has all its orbits
periodic in a neighborhood of the origin, does not depend of the choice of the system of coordinates
on this neighborhood. It is natural to ask if the coefficients of the return mapping are
$S^{1}$-invariant polynomials. The answer is no, but it is possible to find another system of
generators of the Bautin ideal which are $S^{1}$-invariant and so that the associated analytic
series is $\phi$-equivalent to the return mapping. This was proved by ([6]) using the approach of
the successive derivatives.  For $d=2,3$ it was
proved by H. Zoladek ([15])using the classical Bautin approach. We recall here his proof and derive the
result in its full generality for any $d$, $(d \geq 2)$.

        Let us consider the ring of polynomials in
$(a,b), [\theta, ({\rm sin}(\theta),{\rm cos}(\theta)]$:
$$R[a,b]\bigotimes R[\theta, ({\rm sin}(\theta),{\rm cos}(\theta)].$$
 The rotation group acts on the
ring $R[a,b]$ as displayed previously and so it acts on $R[a,b]\bigotimes R[\theta, ({\rm
sin}(\theta),{\rm cos}(\theta)]$ in an obvious way (just extend the action by the identity on the
second factor). We denote $p \mapsto \phi*p$ this action where $\phi$ belongs to $S^{1}$ and $p$
belongs to $R[a,b]\bigotimes R[\theta, ({\rm sin}(\theta),{\rm cos}(\theta)]$. Now $S^{1}$ acts also
on the ring
$$ R[\theta, ({\rm sin}(\theta),{\rm cos}(\theta)]$$

in a natural way:

$$ \phi*(\theta)= \theta+\phi, \phi*({\rm sin}(\theta),{\rm cos}(\theta))=({\rm
sin}(\theta+\phi),{\rm cos}(\theta+\phi)). \eqno(1.34)$$

        This action can be extended to an action $p \mapsto \phi^{*}p$ on the ring
$$R[a,b]\bigotimes R[\theta,
({\rm sin}(\theta),{\rm cos}(\theta)]$$ by the identity on the first factor.
\vskip 12pt
{\bf Definition I.6.1}
\vskip 12pt

        An element $p$ of the ring $R[a,b]\bigotimes R[\theta, ({\rm sin}(\theta),{\rm cos}(\theta)]$ is said
to be covariant if for all $\phi$: $\phi*p = \phi^{*}p$.
\vskip 12pt
{\bf Lemma I.6.2}
\vskip 12pt

        Product and sum of two covariant polynomials are covariant polynomials.
\vskip 12pt

        The proof of this lemma is very easy and is omitted.
\vskip 12pt
{\bf Lemma I.6.3}
\vskip 12pt

        The coefficients $R_{k}(\theta)$ of the equation (1.8) are covariant polynomials.

\vskip 12pt
{\bf Proof.}
\vskip 12pt

        The trigonometric polynomials $A(\theta)$ and $B(\theta)$ are clearly covariant polynomials as
seen from the definition of the actions. Now from Lemma I.6.2 it follows that the products
$A(\theta)B(\theta)^{k}$ are covariant polynomials.

        It is now necessary to be more careful with the inductive construction of the coefficients
$v_{k}(\theta)$. We initiate the discussion with the term of lowest degree $v_{d}(\theta)$ defined
by the condition:

$$v_{d}'(\theta) = R_{d}(\theta). \eqno(1.35)$$

        The term $R_{d}(\theta)$ to be integrated is a trigonometric polynomial. Thus this yields:

$$v_{d}(\theta) = z_{d}(\theta) + s_{d}(\theta), \eqno(1.36)$$

{\noindent}where $z_{d}$ is a constant and $s_{d}(\theta)$ is a trigonometric polynomial with
$s_{d}(0)=s_{d}(2\pi)=0$. Note that

$$z_{d} = (1/{2\pi})\int_{0}^{2\pi}{R_{d}(\phi)d\phi}, \eqno(1.37)$$

and the first coefficient of the return mapping is:

$$v_{d}(2\pi) = 2\pi z_{d} =\int_{0}^{2\pi}{R_{d}(\phi)d\phi}. \eqno(1.38)$$

        Next step of the recurrence:

$$ v_{d+1}'(\theta) = \sum^{d+1}_{i=d} B_{i d+1}[v_{d}(\theta)]R_{i}(\theta), \eqno(1.39)$$

        leads to:

$$v_{d+1}(\theta) = z_{d+1} \theta + r^{(2)}_{d+1} {\theta}^{2} + s_{d+1}(\theta), \eqno(1.40)$$

where $s_{d+1}(\theta)$ is a trigonometric polynomial so that: $s_{d+1}(0)=s_{d+1}(2\pi) = 0$,

$$ z_{d+1} = (1/2\pi)\int_{0}^{2\pi}[ \sum^{d+1}_{i=d} B_{i d+1}[s_{d}(\phi)]R_{i}(\phi)]d\phi,
\eqno(1.41)$$

and where the coefficient $R_{d+1}$ is proportional to the coefficient $v_{d+1}$.

        The general step $k$ of the recurrence displays:

$$v_{k}(\theta) =  z_{k} \theta + r^{(2)}_{k} {\theta}^{2}+...+ r^{(k-d+1)}_{k} {\theta}^{k-d+1} +
s_{k}(\theta), \eqno(1.42)$$

where  $s_{k}(\theta)$ is a trigonometric polynomial so that: $s_{k}(0)=s_{k}(2\pi) = 0$,

such that:

$$ s_{k}'(\theta) = \sum^{k}_{i=2} B_{ik}[s_{d}(\theta), ... , s_{k-1}(\theta)]R_{i}(\theta),
\eqno(1.43)$$

$$ z_{k} = (1/2\pi)\int_{0}^{2\pi}[\sum^{k}_{i=2} B_{ik}[s_{d}(\phi), ... ,
s_{k-1}(\phi)]R_{i}(\phi)]d\phi,\eqno(1.44)$$

and where the coefficients $ r^{(j)}_{k}, j=2,..., k-d+1$ belongs to the ideal generated by the
preceding coefficients $z_{d},...,z_{k-1}$. The construction of the coefficients $z_{k}$ yields:
\vskip 12pt
{\bf Lemma I.6.4}
\vskip 12pt

        The two series $L(x)-x = \sum _{k \geq d} v_{k}(2\pi)x^{k}$ and $\Phi(x)-x = \sum _{k \geq d}
z_{k}x^{k}$ are $\phi$-equivalent.
\vskip 12pt

        We prove now that the polynomials $v_{k}$ are $S^{1}$-invariant. The proof splits into several
lemmas of independent interest.
\vskip 12pt
{\bf Lemma I.6.5}
\vskip 12pt

        Let $T_{k}(\theta)$ be an element of the ring $R[a,b]\bigotimes R[\theta, ({\rm sin}(\theta),{\rm
cos}(\theta)]$ which is a covariant trigonometric polynomial, then the polynomial $v_{k}$
defined by:

$$v_{k} = \int_{0}^{2\pi}T_{k}(\phi)d\phi, \eqno(1.45)$$

is an invariant polynomial.
\vskip 12pt
{\bf Proof.}
\vskip 12pt

        Let $\theta$ be an element of $S^{1}$, consider $\theta*(v_{k})$:

$$\theta*(v_{k}) = \int_{0}^{2\pi}{\theta*T_{k}}(\phi)d\phi. \eqno(1.46)$$

        The covariance property of $T_{k}$ yields:

$$ \int_{0}^{2\pi}{\theta*T_{k}}(\phi)d\phi= \int_{0}^{2\pi}{\theta^{*}T_{k}}(\phi)d\phi=
 \int_{0}^{2\pi}{T_{k}}(\phi+\theta)d\phi =  \int_{0}^{2\pi}{T_{k}}(\phi)d\phi. \eqno(1.47)$$

\vskip 12pt
{\bf Lemma I.6.6}
\vskip 12pt

        Let $T_{k}(\theta)$ be a trigonometric polynomial without constant term which is covariant.
 The only trigonometric polynomial $S_{k}(\theta)$ without constant term so that:

$$S_{k}'(\theta) = T_{k}(\theta), \eqno(1.48)$$

is covariant.
\vskip 12pt
{\bf Proof.}
\vskip 12pt

        Write

$$T_{k}(\theta) = \sum_{l} T_{kl}(a,b)e^{il\theta}. \eqno(1.49)$$

        Let us consider the first action of $S^{1}$. Let $\phi$ be an element of $S^{1}$. Under the action
 of this element, each term $T_{kl}(a,b)$ gets multiplied by the factor $e^{i\phi}$. Now the
polynomial $S_{k}(\theta)$ writes:

$$S_{k}(\theta) = \sum_{l} [T_{kl}(a,b)/l]e^{il\theta}. \eqno(1.50)$$

        Its covariance follows easily from this expression.
\vskip 12pt

 We conclude this paragraph with the theorem
\vskip 12pt
{\bf Theorem I.6.7}
\vskip 12pt
        The displacement function $L(x)-x = \sum _{k \geq d} v_{k}(2\pi)x^{k}$ is $\phi$-equivalent to a
series $\Phi(x)-x = \sum _{k \geq d} z_{k}x^{k}$ with $S^{1}$-invariant coefficients $z_{k}$.

\vskip 12pt
{\bf Proof.}
\vskip 12pt

        The polynomial $R_{d}(\theta)$ is covariant (cf. Lemma I.6.3) and thus the coefficient $z_{d}$ is
invariant as a consequence of formula (1.44) and of the Lemma I.6.5. The polynomial $s_{d}$ is
covariant as a result of the Lemma I.6.6. Assume inductively that the polynomials
$s_{d}{\theta},...,s_{k-1}(\theta)$ are covariant. Then  Lemma I.6.3 and Lemma I.6.6 and the
formula (1.43) yield the fact that $s_{k}$ is covariant and the Lemma I.6.5 imply that the
coefficient $z_{k}$ is invariant.

\vskip 12pt
{\bf II-A GLOBAL AND COMPLEX GENERALIZATION OF BAUTIN'S THEOREM.}
{\vskip 12pt}

        Bautin's method as displayed in the preceding chapter does not allow to produce a precise estimate
on the size of the domain on which the number of zeros of an $A_{0}$-series (such as the
displacement function) is controlled. It is not surprising that to get such an estimate, complex
analytic methods are more appropriated. We begin this chapter with general considerations on
Bernstein classes.
\vskip 12pt
\noindent

{\bf II.1 Finiteness properties, Bernstein inequality.}
{\vskip 12pt}

        In this part, the definitions of the Bernstein classes $B_{1}$
and $B_{2}$ are displayed and their equivalence is proved. Then it
is shown that the number of zeros of functions of these classes is
explicitly bounded. The closed disc centered at $0\in C$ of radius
$R$ is denoted ${\bar D}_{R}$. This paragraph relies entirely on
([8]) and ([9]).
{\vskip 12pt}

{\bf Definition II.1.1}
{\vskip 12pt}

        Let $R>0$, $0<\alpha<1$ and $K>0$ and $f$ holomorphic in a
neighborhood of ${\bar D}_{R}$. The function $f$ belongs to the
Bernstein class $B_{R, \alpha, K}^{1}$ if and only if:
$${\rm max}[{\mid f(z)\mid}, z\in{\bar D}_{R}]/{\rm max}[{\mid f(z)\mid}
, z\in{\bar D}_{\alpha R}]
\leq K. \eqno(2.1)$$

{\bf Definition II.1.2}
{\vskip 12pt}

        Let $N$ be an integer, $R>0$ and $c>0$ and $f(z)=\sum_{i\geq 0} a_{i}z^{i}$
be an analytic function on a neighborhood of $0\in C$. The
function $f$ belongs to the Bernstein class $B_{N, R, c}^{2}$ if
and only if for all $j\geq N$:

$$\mid a_{j}\mid R^{j}\leq c {\rm max}_{i=0,...,N}(\mid a_{i}\mid R^{i}). \eqno(2.2)$$

{\bf Theorem II.1.3}
{\vskip 12pt}

        Let $f$ be an element of $B_{N, R, c}^{2}$. Then $f$ is
analytic on the open disc $D_{R}$, and for all $R'<R$ and
$\alpha<1$, and $K={\alpha^{-N}}[1+\alpha(1-{\alpha}^{N})/(1-\alpha)+c\beta/(1-\beta)]$
where $\beta=R'/R$, $f$ belongs to $B_{R', \alpha, K}^{1}$.
        Conversely if $f$ belongs to $B_{R, \alpha, c}^{1}$, then $f\in B_{N, R, c}^{2}$
with $N= [({\rm log}_{2}(K)-{\rm log}_{2}(1-\alpha)+1)/({\rm
log}_{2}(1/\alpha))]$, $c=[K(2K+1)/(1-\alpha)^{2}$.
{\vskip 12pt}

{\bf Proof}
{\vskip 12pt}

        If $f\in B_{N, R, c}^{2}$, the convergence of the series $f(z)=\sum_{i\geq 0} a_{i}z^{i}$
on the disc $D_{R}$ is consequence of the Cauchy-Hadamard formula
and of the inequality:
$$\mid a_{n}\mid \leq c{\rm max}_{i=0,...,N}(\mid a_{i}\mid R^{i})(R^{-n}),\eqno(2.3)$$
for large $n$.
        Denote $m={\max}(\mid f(z)\mid, z\in D_{\alpha R}$,
Cauchy's inequalities yield:
$$\mid a_{i}\mid\leq m/(\alpha R')^{i}.\eqno(2.4)$$ For $i=0,...,N$, this
displays:
$$\mid a_{i}\mid R^{i}\leq m/(\alpha R'/R)^{i}\leq [m/{\alpha}^{N}(R'/R)^{N}]\eqno(2.5)$$
and then for all $j \geq N+1$:
$$\mid a_{j}\mid R^{j}\leq c[m/{\alpha}^{N}(R'/R)^{N}].\eqno(2.6)$$
This yields:
$${\max}(\mid f(z)\mid, z\in D_{R'}\leq \sum_{i=0}^{N}|a_{i}|R'^{i}+
\sum_{j\geq N+1}|a_{j}|R'^{j}\eqno(2.7)$$
$$\leq m\sum_{i=0}^{N}(1/\alpha R')^{i}R'^{i}
+c[m/{\alpha}^{N}(R'/R)^{N}]\sum_{j\geq N+1}(R'/R)^{j}\eqno(2.8)$$,
and thus:
$${\max}(\mid f(z)\mid, z\in D_{R'})\leq [m/(\alpha^{N})]
[1+\alpha(1-{\alpha}^{N})/(1-\alpha)+c\beta/(1-\beta)]\eqno(2.9)$$
as expected.
        Conversely, assume that $f$ belongs to $B_{R, \alpha,
K}^{1}$, $f(z)=\sum_{i\geq 0} a_{i}z^{i}$. Let $N$ be an integer
and $\sigma={\rm max}_{i=0,...,N}(|a_{i}|R^{i})$ and
$b={\rm max}(|R(z)|, z\in D_{\alpha R})$, where $R(z)=\sum_{i\geq
N+1}a_{i}z^{i}$. Consider $P_{N}(z)=\sum_{i=0}^{N}a_{i}z^{i}$ and:
$${\rm max}(|P_{N}(z)|, z\in D_{\alpha R})\leq \sum_{i=0}^{N}|a_{i}|(\alpha R)^{i}\eqno(2.10)$$
$${\rm max}(|P_{N}(z)|, z\in D_{\alpha R})\leq
\sigma/(1-\alpha).\eqno(2.11)$$
This yields:
$${\rm max}(|f(z)|, z\in D_{\alpha R})\leq
K[b+(\sigma/(1+\alpha)].\eqno(2.12)$$
Cauchy inequalities yield:
$$|a_{j}|R^{j}\leq K[b+(\sigma/(1+\alpha)].\eqno(2.13)$$
The following bound for $b$ may be displayed:
$$b\leq \sum_{j \geq N+1}[|a_{j}|(\alpha R)^{j}\leq K[b+(\sigma/(1+\alpha)][\alpha^{N+1}
/(1-\alpha)],\eqno(2.14)$$
$$b[1-(K\alpha^{N+1}/(1-\alpha)]\leq
K\sigma[\alpha^{N+1}/(1-\alpha)^{2}].\eqno(2.15)$$
Fix $N$ so that:
$$(K\alpha^{N+1}/(1-\alpha)\leq 1/2;\eqno(2.16)$$
this is achieved for instance if:
$N= [({\rm log}_{2}(K)-{\rm log}_{2}(1-\alpha)+1)/({\rm
log}_{2}(1/\alpha))]$. Such a choice for $N$ yields:
$$b\leq 2K\sigma[\alpha^{N+1}/(1-\alpha)^{2}]\leq 2K\sigma/(1-\alpha)^{2}.\eqno(2.17)$$
Now this yields:
$$|a_{j}|R^{j}\leq K{{[2\sigma K/(1-\alpha)^{2}]+[\sigma/(1-\alpha)]}},\eqno(2.18)$$
$$|a_{j}|R^{j}\leq K(2K+1)/(1-\alpha)^{2};\eqno(2.19)$$
the theorem now follows with:
$c=[K(2K+1)/(1-\alpha)^{2}$.

{\vskip 12pt}

        The number of zeros of a function which belongs to a
Bernstein class is bounded accordingly to what follows.

{\vskip 12pt}
{\bf Theorem II.1.4}
{\vskip 12pt}
        If $f\in B_{R, \alpha, K}^{1}$, then the number of zeros
of $f$ in the disc ${\bar D}_{\alpha R}$ is less than:
$${\rm log}_{2}(K)/[{\rm log}_{2}[(1+\alpha^{2})/2\alpha]]\eqno(2.20)$$
{\bf Proof:}
{\vskip 12pt}
        Assume that $f$ displays $n$ zeros $z_{1},..., z_{n}$ in the disc
${\bar D}_{\alpha R}$. Denote $g$ the holomorphic function on the
disc $D_{R}$ defined as:
$$g(z)= f(z)\prod_{k=1}^{n}[(R^{2}-z{\bar z}_{k})/R(z-z_{k})].\eqno(2.21)$$
The maximum principle yields:
$${\rm max}(|g(z)|, z\in D_{R})={\rm max}(\mid f(R{\rm e}^{i\theta})
\prod_{k=1}^{n}[(R^{2}-R{\rm e}^{i\theta}{\bar z}_{k})
/R(R{\rm e}^{i\theta}-z_{k})]\mid, \theta\in [0, 2\pi]),\eqno(2.22) $$
$$= {\rm max}(\mid f(R{\rm e}^{i\theta})\mid \mid
\prod_{k=1}^{n}[(R^{2}-R{\rm e}^{i\theta}{\bar z}_{k})
/R(R-{\rm e}^{i\theta}{\bar z}_{k})]\mid, \theta\in [0, 2\pi]),\eqno(2.23)$$
$$={\rm max}(|f(z)|, z\in D_{R}).\eqno(2.24)$$
Consider now:
$${\rm max}(|g(z)|, z\in D_{\alpha R})={\rm max}(\mid f(\alpha R{\rm
e}^{i\theta}))\prod_{k=1}^{n}[(R^{2}-\alpha R{\rm e}^{i\theta}{\bar z}_{k})
/R(\alpha R{\rm e}^{i\theta}-z_{k})]\mid) \eqno(2.25)$$

$$\geq {\rm max}(|f(z)|, z\in D_{\alpha R})
\prod_{k=1}^{n}{\rm min}[\mid R^{2}-\alpha R{\rm e}^{i\theta}{\bar z}_{k})
/R(\alpha R{\rm e}^{i\theta}-z_{k})\mid.\eqno(2.26)$$
        Write now separately each quantities:

$$\mid (R^{2}-\alpha R{\rm e}^{i\theta}{\bar z}_{k})
/R(\alpha R{\rm e}^{i\theta}-z_{k})]\mid\eqno(2.27)$$
and $z_{k}=R_{k}{\rm e}^{i\theta_{k}}$. The minimum of the
quantity when $\theta-\theta_{k}$ varies in $[0, 2\pi]$ is
$$(R+\alpha R_{k})/(\alpha R+R_{k}).\eqno(2.28)$$
The minimum of this quantity when $0\leq R_{k}\leq \alpha R$ is
$$(1+{\alpha}^{2})/2\alpha.\eqno(2.29)$$
        This yields the following inequality:
$${\rm max}(|g(z)|, z\in D_{\alpha R})\geq {\rm max}(|f(z)|, z\in D_{\alpha R})
[(1+{\alpha}^{2})/2\alpha]^{n}.\eqno(2.30)$$
        The inequalities:
$${\rm max}(|f(z)|, z\in D_{R})={\rm max}(|g(z)|, z\in D_{R})\geq
{\rm max}(|g(z)|, z\in D_{\alpha R}),\eqno(2.31)$$
display:
$$K\geq [(1+{\alpha}^{2})/2\alpha]^{n},\eqno(2.32)$$
and thus:
$$n\leq {\rm log}_{2}(K)/[{\rm log}_{2}[(1+\alpha^{2})/2\alpha]].\eqno(2.33)$$
\vskip 12pt
{\bf Lemma II.1.5}
\vskip 12pt
        Let $f\in B_{N, R, c}^{2}$ and $R''=R/[2^{3N}{\rm max}(c,
2)]$, the number of zeros of $f$ in the disc ${\bar D}_{R''}$ is
less than $N$.
{\vskip 12pt}
{\bf Proof}.
{\vskip 12pt}
        Set $R'=R/[{\rm max}(c,2)]$ and $\alpha=2^{-3N}$. This
displays:
$$[(1+{\alpha}^{2})/2\alpha]=(2^{6N}+1)/(2^{3N+1})>2^{3N-1},\eqno(2.34)$$
$$K<2^{3N.N+2}.$$
This yields that the number of zeros of $f$ is less than:
$$(3N^{2}+2)/(3N-1)<N+1,\eqno(2.35)$$
as soon as $N\geq 2.$

\vskip 12pt
\noindent
{\bf II.2 The Hironaka polynomial division theorem}
{\vskip 12pt}

        The precise statement of the following theorem was provided by P. Milman. The terminology
"Hironaka polynomial division theorem" is also due to him. This technical result allows to clarify
the presentation of the  "Quantitative Bautin's theorem" first proved in ([8]).
{\vskip 12pt}

{\bf Theorem II.2.1}
{\vskip 12pt}

        Let $I$ be an ideal of $C[\lambda_{1},...,\lambda_{D}]$. There is a system of generators $g_{1},...,
 g_{s}$ of the ideal $I$ and constants $C$ and $C_{1}$ such that for all elements $f$ of degree $k$ of $I$,
 there is a decomposition:

$$f(\lambda) = \sum_{i=1}^{s} \phi_{i}(\lambda)g_{i}(\lambda), \eqno(2.36)$$

with

$$deg(\phi_{i}) \leq deg(f) =k, \eqno(2.37a) $$

and

$$\mid \phi_{i}\mid \leq CC_{1}^{k}\mid f\mid. \eqno(2.37b) $$
{\vskip 12pt}

                We include here a full proof of this theorem. The non classical part of this division theorem
 is displayed in the control of the norms (equation 2.37b) of the quotients.
 The classical proof uses the inversion of an operator defined
 on Banach spaces of analytic functions but we do not follow these lines here.
{\vskip 12pt}
{\bf Definition II.2.2}
{\vskip 12pt}

        In the following, a total ordering $\leq$ on $N^{D}$ is
 said to be compatible with the addition if:

  i) For all indices $\alpha \in N^{D}$, $\beta \in N^{D}$, then $\alpha \leq \alpha+\beta$,
  ii) For all indices $\alpha^{1}, \alpha^{2}, \beta$, $\alpha^{1}+\beta \leq \alpha^{2}+\beta$ if
and only if $\alpha^{1} \leq \alpha^{2}$.
{\vskip 12pt}

 Here, we choose (for instance) the total ordering defined in such manner:

$\alpha \leq \beta$ if and only if $\alpha_{1}+...+\alpha_{D} = \mid\alpha\mid \leq
\beta_{1}+...+\beta_{D} = \mid\beta\mid$ or if $\mid\alpha\mid = \mid\beta\mid$, then
$C(\alpha)\leq C(\beta)$ where $C(\alpha)= \sum_{i=0}^{D}
\alpha_{i}\epsilon^{i}$, $\epsilon$ is a transcendent number, $0\leq \epsilon \leq
1$. Note firstly that $\leq$ defines a total ordering compatible
with the addition. Indeed, the transcendent nature of $\epsilon$
yields
$$C(\alpha)=C(\beta){\Leftrightarrow} \alpha=\beta\eqno(2.38)$$
{\vskip 12pt}

{\bf Definition II.2.3}
{\vskip 12pt}

        Let $f \in C[\lambda], f \neq 0, f = \sum f_{\alpha}\lambda^{\alpha}$. The largest exponent
$\alpha$ so that $f_{\alpha} \neq 0$ is called the privileged exponent of $f$ and is denoted
$exp(f)$. The monomial $[f_{\alpha}\lambda^{\alpha}, \alpha=exp(f)]$ is called the initial monomial
and denoted $In(f)$.
{\vskip 12pt}

{\bf Definition II.2.4}
{\vskip 12pt}

        Given $s$ polynomials $g_{1},...,g_{s}$, the associated partition of $N^{D}$ is defined as follows:

 $\Delta_{1} = exp(g_{1})+N^{D},...,\Delta_{i}= exp(g_{i})+N^{D}-\bigcup_{j<i}\Delta_{j},...,
 \bar\Delta =N^{D}-\bigcup_{i=1}^{s}\Delta_{i}$.
{\vskip 12pt}
{\bf Definition II.2.5}
{\vskip 12pt}
        Let $I$ be an ideal of $C[\lambda_{1},...,\lambda_{D}]$.
 Consider the set:
 $$exp(I)=(exp(f), f{\in}I).\eqno(2.39)$$
 It can be shown that this set has finitely many extremal points:
 $$\alpha^{1},...,\alpha^{s}.\eqno(2.40)$$
 Choose $g_{1},..., g_{s}$ in the ideal $I$ so that $exp(g_{i})=\alpha^{i},
 i=1,...s$. Such a set of polynomials is called a standard basis,
 Hironaka basis or Grobner basis (of the ideal $I$ relatively to
 the ordering $\leq$).
{\vskip 12pt}

 {\bf Proposition II.2.6}
 {\vskip 12pt}

        Let $f \in C[\lambda], deg(f)=k$, and $g_{1},...,g_{s} \in C[\lambda]$, there is a constant $C$ which depends
only on the polynomials $g_{i}$ and there are unique $h_{1},...,h_{s}, h\in C[\lambda]$ such that:

$$ f = h_{1}g_{1}+...+h_{s}g_{s} + h, \eqno(2.41a)$$
$$ h_{i}=\sum h_{i,\alpha}\lambda^{\alpha} \Rightarrow \alpha+exp(g_{i})\in\Delta_{i}, \eqno(2.41b)$$
$$h =\sum h_{\alpha}\lambda^{\alpha} \Rightarrow \alpha\in\bar\Delta, \eqno(2.41c)$$
$$Max(\mid h_{i}\mid,\mid h \mid) \leq CC_{1}^{k}\mid f \mid,\eqno(2.41d)$$
$$Max[deg(h_{i}), deg(h)] \leq k.\eqno(2.41e)$$
{\vskip 12pt}

{\bf Proof:}
{\vskip 12pt}

        Given $f$, denote $In(f)=f_{\alpha_{0}}{\lambda}^{\alpha_{0}}$.
If $\alpha_{0}\in \bar\Delta$ set $h^{(1)}=f_{\alpha_{0}}{\lambda}^{\alpha_{0}}$
and $h_{i}^{(1)}=0, i=1,...,s$, $f^{(1)}=f-h^{(1)}$. This yields $\mid f^{(1)}
\mid\leq\mid f \mid$
, $deg(f^{(1)})\leq deg(f)$ and $\mid h^{(1)}\mid\leq\mid f \mid$
, $deg(h^{(1)})\leq deg(f)$.
If $\alpha_{0}\in\Delta_{i}$ (it may belong to several $\Delta_{j}$, choose
one). Then set:
$h^{(1)}_{j}=0$ if $j\neq i$ and $h^{(1)}=0$. Write
$\alpha_{0}=exp(g_{i})+\beta_{0}$,
$In(g_{i})=\gamma_{i}{\lambda}^{exp(g_{i})}$, then denote
$h_{i}^{(1)}=f_{\alpha_{0}}{\lambda}^{\beta_{0}}/{\gamma_{i}}$ and
$f^{(1)}=f-h_{i}^{(1)}g_{i}$. This yields:
$$\mid h^{(1)}_{i}\mid \leq C\mid f \mid\eqno(2.42)$$
where:
$$C={\rm max}(1/\gamma_{j}), j=1,...,s\eqno(2.43)$$
$${\rm deg}(h^{(1)}_{i}) = \mid \beta_{0}\mid\leq k-{\rm min}(\mid exp(g_{j})\mid)\eqno(2.44)$$
and thus this displays:
$$\mid f^{(1)}\mid\leq\mid f \mid+CG\mid f \mid\leq (1+CG)\mid f \mid\eqno(2.45)$$
where:
$$G={\rm max}(\mid g_{i} \mid), i=1,...,s.$$ Note furthermore
that:
$${\rm deg}(f^{(1)})\leq {\rm max}[{\rm deg}(f), {\rm deg}(h_{i}^{(1)}g_{i})]\eqno(2.46)$$
and
$${\rm deg}(h_{i}^{(1)}g_{i})= \mid \beta_{0} \mid+{\rm deg}(g_{i})\leq
\mid \beta_{0} \mid+\mid exp(g_{i})\mid\leq\mid \alpha_{0}\mid\leq k\eqno(2.47)$$
and thus
$${\rm deg}(f^{(1)})\leq k\eqno(2.48).$$
        Then, repeat the whole process with $f^{(1)}$ in place of $f$.
Note that the privileged exponent of $f^{(1)}$ is strictly less
than the privileged exponent of $f$. The process stops
ultimately when the privileged exponent becomes zero. The
choice of the ordering yields that the number of steps involved is
less than $k(1+{\epsilon}^{1-D})$. This yields the result of the
proposition with:
$$C_{1}=(1+CG)^{1+{\epsilon}^{1-D}}\eqno(2.49)$$.

        The theorem follows of the preceding proposition and of
the property that if a polynomial $f$ belongs to an ideal $I$ and
if $g_{1},..., g_{s}$ is a Grobner basis of the ideal $I$, then
the division of $f$ yields $h=0$. (See for instance
Lejeune-Jalabert ([10])).
\vskip 12pt
\noindent

{\bf II.3 The main theorem}
{\vskip 12pt}

{\bf Theorem II.3.1}
{\vskip 12pt}

        Let $f_{\lambda}(x) = \sum_{k \geq 1} f_{k}(\lambda)x^{k}$ be an $A_{0}$-series.
 For any $\lambda \in C^{D}$, the function $f_{\lambda}$ belongs to the Bernstein class
 $B_{d-1,R,c}^{2}$ where $d$ is the Bautin index of the series $f_{\lambda}(x)$,
 $R= [(C_{1}{\bar{\lambda}})^{K_{1}}K_{4}]^{-1}$, $c=
[MCK_{3}(C_{1}{\bar{\lambda}})^{K_{2}}]/R^{d}$, if $R \leq 1$, $c=
[MCK_{3}(C_{1}{\bar{\lambda}})^{K_{2}}]/R$, if $R \geq 1$,
 $\bar{\lambda}= max(1, \mid\lambda \mid)$. The constants $K_{1},...,K_{4}$ are those which appear in
the definition (I.2.1) of an $A_{0}$-series. The constants $C$ and $C_{1}$ appears in the "Hironaka
polynomial
division theorem" and $M$ is the norm of the matrix  of the change of basis $(f_{k})$ to
a Grobner basis $(g_{l})$ of the Bautin ideal of the series $f_{\lambda}(x)$.
{\vskip 12pt}
{\bf Proof}
{\vskip 12pt}
        Write first the condition for $f_{\lambda}(x)$ to be an
$A_{0}$-series as follows:
$$deg[f_{k}(\lambda)]\leq K_{1}k+K_{2},$$
$$\mid f_{k}\mid \leq K_{3}K_{4}^{k}.$$
        Denote $I$ the ideal generated by the $d$ first
coefficients $f_{1}(\lambda),..., f_{d}(\lambda)$ and write:
$$f_{k}(\lambda)=\sum_{j=1}^{s}[\phi_{j}^{k}(\lambda)g_{j}(\lambda)]. \eqno(2.50)$$
        The Hironaka polynomial division theorem yields:
$$\mid \phi_{j}^{k}(\lambda)\mid \leq CC_{1}^{(K_{1}k+K_{2})}\mid
f_{k}\mid{\bar \lambda}^{(K_{1}k+K_{2})}
\eqno(2.51a)$$
$$\leq [CK_{3}C_{1}^{K_{2}}][C_{1}^{K_{1}}K_{4}]^{k}{\bar \lambda}^{(K_{1}k+K_{2})}
. \eqno(2.51b)$$
        The relation:
$$g_{j}(\lambda)=\sum_{l=1}^{d} \Phi_{jl}(\lambda)f_{l}(\lambda), \eqno(2.52)$$
yields:
$$\mid f_{k}(\lambda)\mid
\leq  [CK_{3}C_{1}^{K_{2}}{\bar \lambda}^{K_{2}}][C_{1}^{K_{1}}K_{4}
{\bar \lambda}^{K_{1}}]^{k}
\sum_{j=1}^{s}\mid g_{j}(\lambda)\mid, \eqno(2.52)$$
$$\leq M[CK_{3}C_{1}^{K_{2}}{\bar \lambda}^{K_{2}}][C_{1}^{K_{1}}K_{4}
{\bar \lambda}^{K_{1}}]^{k}{\rm max}
[\mid f_{1}(\lambda)\mid ,..., \mid f_{d}(\lambda)\mid]. \eqno(2.53)$$
        Write now:
$$R= [(C_{1}{\bar \lambda})^{K_{1}}K_{4}]^{-1}, \eqno(2.54)$$
the equation (2.53) yields:
$$\mid f_{k}(\lambda)\mid R^{k}\leq [MCK_{3}(C_{1}{\bar
\lambda})^{K_{2}}/R^{d}]
{\rm max}
[\mid f_{1}(\lambda)\mid ,..., \mid f_{d}(\lambda)\mid]R^{d}. \eqno(2.55)$$
If $R\leq 1$, this displays:
$$\mid f_{k}(\lambda)\mid R^{k}\leq [MCK_{3}(C_{1}{\bar
\lambda})^{K_{2}}/R^{d}]
{\rm max}
[\mid f_{1}(\lambda)\mid R,..., \mid f_{d}(\lambda)\mid R^{d}]. \eqno(2.56)$$
If $R \geq 1$, the equation (2.53) yields:
$$\mid f_{k}(\lambda)\mid R^{k}\leq [MCK_{3}(C_{1}{\bar
\lambda})^{K_{2}}/R]
{\rm max}
[\mid f_{1}(\lambda)\mid ,..., \mid f_{d}(\lambda)\mid]R, \eqno(2.57)$$
and thus:
$$\mid f_{k}(\lambda)\mid R^{k}\leq [MCK_{3}(C_{1}{\bar
\lambda})^{K_{2}}/R]
{\rm max}
[\mid f_{1}(\lambda)\mid R,..., \mid f_{d}(\lambda)\mid R^{d}]. \eqno(2.58)$$
        This concludes the proof.
\vskip 12pt
        The theorem II.3.1 and the lemma II.1.5 now implies the
following:
\vskip 12pt
{\bf Theorem II.3.2}
\vskip 12pt
        The number of zeros of $f_{\lambda}(x)$ in the disc ${\bar D}_{R''}$
is less than $d-1$ with
$$R''=R/[2^{3(d-1)}{\rm max}(c,2)]. \eqno(2.59)$$

\vskip 12pt
{\bf III. Extension to any dimension of Bautin's theory}
\vskip 12pt
   These last years, the dynamics of plane systems was extensively
studied and several new techniques were developed. Some are
specific to $2$-dimensional systems but often these methods
can be appropriately extended to multi-dimensional systems.
        The algorithm of the successive derivatives was derived
 some years ago ([5]) to find the first
non-vanishing derivative (relatively to the parameter $\epsilon$)
of the return mapping (near the origin) of a plane vector field
$X_{0}+{\epsilon}X_{1}$ of type:

$$X_{0}+{\epsilon}X_{1}=x{\partial}/{\partial y}-y{\partial}/{\partial x}+
{\epsilon} \sum^{d} _{i,j/ i+j=2}
 [ a_{i,j} x^{i}y^{j}{\partial}/{\partial x} +
  b_{i,j} x^{i}y^{j}{\partial}/{\partial y}].
  \eqno(3.1)$$

        The algorithm was then used in the center-focus problem
(cf. [6]), which directly relates to Hopf
bifurcations of higher order and to several other problems on
limit cycles of plane vector fields.
        How to extend appropriately this situation in any
 dimension? We have to perturb a dynamics which is integrable
 and displays only periodic orbits. Assuming that the
 perturbation depends of finitely many parameters (say for
 instance it is polynomial), we expect also that the perturbed
 system displays a first return-mapping which is analytic with
 a Taylor expansion with coefficients which depend polynomially of
 the parameters. This return-mapping should label all the periodic
 orbits (at first return) of the perturbed system by its fixed
 points. The principal aim of this paragraph is to present a
 framework where such demands are realized. In this framework, a
 generalization of the algorithm of the successive derivatives is
 provided.
 \bigskip

{\bf III-1 Controlled Nambu dynamics and (*)-property.}
\bigskip
         Let $f=(f_{1},...,f_{n-1}): R^{n} \rightarrow R^{n-1}$ be a
generic submersion (meaning that $f$ is a submersion outside a critical set $f^{-1}(C)$,
where $C$ is a set of isolated points). Let $\Omega =dx_{1}{\wedge}dx_{2}{\wedge}...dx_{n}$ be
a volume form on $R^{n}$. Consider the vector field $X_{0}$ such
that:
$${\iota}_{X_{0}}dx_{1}{\wedge}dx_{2}{\wedge}...{\wedge}dx_{n}=
df_{1}{\wedge}...{\wedge}df_{n-1}.     \eqno (3.2)$$

        The functions $f_{i}, (i=1,...,n-1)$ are first integrals
of the vector field $X_{0}$:

$$df_{i}{\wedge}{\iota}_{X_{0}}dx_{1}{\wedge}dx_{2}{\wedge}...{\wedge}dx_{n}=
(X_{0}.f_{i})dx_{1}{\wedge}dx_{2}{\wedge}...{\wedge}dx_{n}=
df_{i}{\wedge}df_{1}{\wedge}df_{2}{\wedge}...{\wedge}df_{n-1}=0. \eqno(3.3)$$

        This type of dynamics is well-known in Physics and named
Nambu's dynamics ([12]).

        For $c$ varying in a neighborhood of $0$,
assume that the curves $f^{-1}(c)$ have a compact connected
component $\gamma_{c}$. Let $\Sigma$ be a small neighborhood of
the zero-section of the normal bundle to $\gamma_{0}$. For $c$
small enough, the curves $\gamma_{c}$ are closed periodic orbits
of $X_{0}$ and they cut transversely $\Sigma$. Choose $c$ as a
coordinate on the transverse section $\Sigma$ to the flow of
$X_{0}$.
        Lastly, assume that there are $2$-forms $\omega_{i}$ such
that:

$$\iota_{X_{0}}\omega_{i}=df_{i}; i=1,...,n-1. \eqno(3.4)$$

        There are of course different possible choice of the forms
$\omega_{i}$ and accordingly different possible perturbations have to be considered.
If the condition (3.4) is fulfilled, we will
say that the singularity of the Nambu dynamics (3.3) is controlled (or alternatively
that the Nambu dynamics itself is controlled).

        The appropriated extension of the (*)-property first
discussed in ([5]) is presented in the following.

\bigskip
{\bf Definition III.1.1}
\bigskip

        Let $f=(f_{1},...,f_{n-1}): R^{n} \rightarrow R^{n-1}$ be
a generic submersion. Assume that $f^{-1}(c)$ contains a compact
curve $\gamma_{c}$. The application displays the (*)-property if
for all polynomial $1$-forms $\omega$ such that

$$\int_{{\gamma}_{c}}{\omega} = 0,  \eqno(3.5)$$

for all $c$; there exist polynomial $g_{i}, R$ such that:

$$\omega= g_{1}df_{1}+...+g_{n-1}df_{n-1}+dR. \eqno(3.6)$$

It was proved in ([5]) that the function $f_{1}:R^{2} \rightarrow R^{1},
f_{1}: (x_{1}, x_{2}) \rightarrow (x_{1}^{2}+x_{2}^{2})$ displays
the (*)-property. Several generalizations were proposed after but
the core of the argument in the computation of the successive
derivatives is captured in this notion.
        The generalization proposed in this article provides a new
presentation of the (*)-property which seems interesting as well
for the $2$-dimensional case. Indeed, the definition of the vector
field $X_{0}$ given in the preceding introduction yields the:
\bigskip
{\bf Proposition III.1.2}
\bigskip

        Let $\omega$ be a $1$-form such that $\omega(X_{0})=0$,
then there are functions $g_{1},...,g_{n-1}$ so that:

$$\omega = g_{1}df_{1}+...+g_{n-1}df_{n-1}. \eqno(3.7)$$

        Note that the condition $\omega(X_{0})=0$, equivalent to
$\omega{\wedge}df_{1}{\wedge}...{\wedge}df_{n-1}=0$, yields
$\omega = g_{1}df_{1}+...+g_{n-1}df_{n-1}$ where the coefficients $g_{k}$
are obtained as ratio of minors of the Jacobian matrix of the
$f_{j}$.

        This displays an alternative to the (*)-property now
presented as follows:
\bigskip

{\bf Proposition III.1.3}
\bigskip

        A generic submersion $f:R^{n}\rightarrow R^{n-1}$ displays
the (*)-property if for any polynomial $1$-form $\omega$ such that

$$\int_{{\gamma}_{c}}{\omega} = 0,  \eqno(3.5)$$

for all $c$; then there exists a polynomial $R$ such that:

$$\omega(X_{0})=X_{0}.R. \eqno(3.8)$$

        Such a function $R$ can be (in principle) constructed with
the following pattern. Choose $R$ arbitrarily on the transverse
section $\Sigma$, then extend $R$ to the whole tubular
neighborhood of $\gamma_{0}$ saturated by the orbits $\gamma_{c}$
by integration of the $1$-form $\omega$ along the orbits of
$X_{0}$.

\bigskip

{\bf III-2 The successive derivatives of the first return mapping
 of the perturbed system.}

\bigskip

        Now perturb $X_{0}$ into $X_{\epsilon}=X_{0}+{\epsilon}X_{1}$.
Let $M$ be a point of ${\Sigma}$ close to $0$ and let $\gamma_{\epsilon}$
be the trajectory of $X_{\epsilon}$ passing by the point $M$. The
next first intersection point of $\gamma_{\epsilon}$ with $\Sigma$
defines the so-called first return mapping of $X_{\epsilon}$
relatively to the transverse section $\Sigma$: $c \longmapsto L(c,
\epsilon)$.
         The mapping $L$ is analytic and it
displays a Taylor development (in $\epsilon$):

$$L(c,\epsilon)= c+{\epsilon}L_{1}(c)+...+{\epsilon}^{k}L_{k}(c)+O(\epsilon)^{k+1}
. \eqno(3.9)$$

        The expression of the first coefficient $L_{1}(c)$ is
classical and belongs to the lore of bifurcation theory. With the
vector field $X_{\epsilon}$ and the $1$-forms $\omega_{i}$ (cf.
[F]), introduce the $1$-forms:

$$\iota_{X_{\epsilon}}{\omega_{i}}=\iota_{X_{0}}{\omega_{i}}+{\epsilon}
\iota_{X_{1}}{\omega_{i}}=df_{i}+{\epsilon}\iota_{X_{1}}{\omega_{i}}
. \eqno(3.10)$$
\bigskip
{\bf Definition III.2.1}
\bigskip
        The perturbation $X_{\epsilon}$ of the controlled Nambu
dynamics is said to be admissible if for all the 2-forms
$\omega_{i}$, the 1-forms $\iota_{X_{1}}{\omega_{i}}$ have polynomial coefficients.
Note that has said above this displays different admissible
perturbations depending of the choice of the forms $\omega_{i}$.
\bigskip
        Recall that the parameter $c$ chosen as coordinates on
the transverse section $\Sigma$ is the restriction of the
functions $f=(f_{1},..., f_{n-1})$ to the section.

        Then the $i^{th}$-component of $L_{1}(c)$ is equal to:

$$L_{1, i}(c)= \int_{\gamma_{0}}{\iota_{X_{1}}{\omega_{i}}}. \eqno(3.11)$$

        Assume now that the first derivative $L_{1}(c)$ vanishes
identically and that the submersion $f$ displays the (*)-property
then there exist $g_{ij}$ and $R_{i}$ such that:

$$\iota_{X_{1}}\omega_{i}=\sum_{j} {g_{ij}df_{j}}+dR_{i}. \eqno(3.12)$$

        Following the lines of the algorithm of the successive
derivatives, the expression (3.12) yields:

$$L_{2,i}(c)=-\int_{\gamma_{0}}\sum_{j}{g_{ij}{\iota_{X_{1}}\omega_{j}}}. \eqno(3.13)$$

        This is indeed the second step of a general recursive
scheme which displays as follows:

        Assume that all the $k^{th}$-first derivatives of the
return mapping of the perturbed vector field vanish identically. This
yields:

$$L_{k,i}(c)=\int_{\gamma_{0}}\sum_{j}{g^{k-1}_{ij}}
{\iota_{X_{1}}\omega_{j}}=0. \eqno(3.14)$$

        The (*)-property yields new functions $g^{k}_{ij}, R^{k}$
such that:

$$\sum_{j}{{g^{k-1}_{ij}}\iota_{X_{1}}\omega_{j}}=\sum_{j}{g^{k}_{ij}df_{j}}+dR^{k}
. \eqno(3.15)$$
        This yields the following expression of the
$(k+1)^{th}$-derivative of the return mapping of the perturbation:

$$L_{k+1,i}(c)=\int_{\gamma_{0}}\sum_{j}{{g^{k}_{ij}}
{\iota_{X_{1}}\omega_{j}}}. \eqno(3.16)$$

        The algorithm implies of course the first
\bigskip
{\bf Theorem III.2.2}
\bigskip
        Let $X_{0}$ be a controlled Nambu dynamics which displays
the (*)-property and let $X_{1}$ be an admissible perturbation.
Then the perturbed dynamics $X_{\epsilon}$ has an analytic first
return map. The coefficients of the Taylor expansion of this return
mapping depend polynomially on the coefficients of the
perturbation.

        From the general theory of projections of analytic sets
([11]), it now follows:
\bigskip

{\bf Theorem III.2.3}
\bigskip

        There exists a uniform bound to the number
of isolated periodic orbits, which correspond to fixed points of
the first return mapping of $X_{0}+{\epsilon}X_{1}$ which
intersect the transverse section $\Sigma$ in the neighborhood of
$0$.

        The general framework presented here should of course be
illustrated with specific examples (of dimension larger than $2$).
 Some have been worked out
recently by Seok Hur (Paris VI) and will be matter to further
publications.

\vfill\eject
\REFERENCES

\ref{N.N. Bautin:}{}{On the number of limit cycles which appear with the variation
of coefficients from an equilibrium position of focus or center type.
Amer. Math. Soc. Transl.{\bf 1, 5}, 336-413 (1962). Translated from
Mat. Sbornik, {\bf 30}, 181-196, (1952).}

\ref{M. Briskin, J.-P. Francoise, Y. Yomdin:}{}{The Bautin ideal
of the Abel Equation. Nonlinearity {\bf 11}, 431-443 (1998).}

\ref{M. Briskin, J.-P. Francoise, Y. Yomdin:}{}{Center conditions,
composition of polynomials and moments on algebraic curves.
Ergod. Th.  Dynam. Sys. {\bf 19}, 1201-1220 (1999).}

\ref{M. Briskin, Y. Yomdin:}{}{Algebraic Families of Analytic
Functions. J. Diff. Equations. {\bf 136} (2), 248-267 (1997).}

\ref{J.-P. Francoise:}{}{Successive derivatives of a first-return
map, application to quadratic vector fields. Erg. Th. Dynam.
Sys., {\bf 16}, 87-96 (1996).}

\ref{J.-P. Francoise,R. Pons:}{}{Computer Algebra Methods and
The stability of differential systems. Random and computational
dynamics, {\bf vol. 3, n°4}, 265-287 (1995).}

\ref{J.-P. Francoise, C.C. Pugh:}{}{Keeping track of limit cycles.
J. Differential Equations {\bf 65}, 139-157 (1986).}

\ref{J.-P. Francoise, Y. Yomdin:}{}{ Bernstein inequality and applications
to analytic geometry and differential equations. J. Funct. Analysis {\bf 146},
 185-205 (1997).}

\ref{J.-P. Francoise, Y. Yomdin:}{}{Projection of analytic sets and Bernstein
inequalities. Singularities Symposium-Lojasiewicz 70, Edts B. Jacubczyk,
W. Pawlucki, Y. Stasica, Banach Center Publications, Warszawa {\bf 44},
103-108 (1998).}

\ref{M. Lejeune-Jalabert:}{}{Effectivit\'e de calculs polynomiaux.
Publications de l'Institut Fourier, Universit\'e de Grenoble,
(1984)}

\ref{S. Lojasiewicz, R. Tougeron, M. Zurro.:}{}{Eclatement des
coefficients des series enti\'eres et deux th\'eor\`emes de Gabrielov.
Manuscripta Matematica, {\bf 92}, 325-337, (1997).}

\ref{Nambu}{}{{\it Broken Symmetry}. Selected papers of Y. Nambu
edited and selected with a foreword by T. Eguchi and K. Nishijima.
World scientific series in 20th century physics, {\bf 13}. World scientific
publishing co. inc. (1995)}

\ref{S. Smale:}{}{Dynamics retrospective: great problems, attempts
that failed. Physica D, {\bf 51}, 267-273, (1991).}

\ref{S. Smale:}{}{Mathematical Problems for the Next Century. The
Mathematical Intelligencer, {\bf 20}, 7-15, (1998).}

\ref{Y. Yomdin:}{}{Global finiteness properties of analytic families
and algebra of their Taylor coefficients. The Arnol'd Fest (Toronto,
Ontario, 1997), 527-555, Fields Inst. Commun. {\bf 24},
Amer. Math. Soc., Providence, RI, (1999).}

\ref{H. Zoladek:}{}{Quadratic systems with center and their perturbations.
J. Differential Equations, {\bf 109}, 223-273, (1994).}

\end